\documentclass
{article}[12pt]
\usepackage[cp1251]{inputenc}
\usepackage[russian]{babel}
\usepackage{amsfonts,amssymb,amsmath,euscript,mathrsfs}

\begin{document}
\vskip0.5cm

\begin{center}\textbf{Lower Bound for the Discrete Norm of a Polynomial on the
Circle}\end{center}
 \begin{center}\textbf{V. N. Dubinin}\end{center}

\begin{center}1. INTRODUCTION AND STATEMENT OF THE RESULT\end{center}
For the approximation of functions, a uniform grid of values of
the arguments is often chosen. In this connection, it is natural
to pose the question of how the discrete norm on a given grid
relates to the uniform norm of the corresponding function on a
given set. In the comparatively recent paper [1], Sheil- Small
showed that, for algebraic polynomials $P$ of degree $n$ and
natural $N>n$, the following estimate holds:
$$
\max\limits_{\omega^N=1}|P(\omega)|\geq\sqrt{\frac{N-n}{n}}
\max\limits_{|z|=1}|P(z)|.\eqno(1)
$$

Earlier Rakhmanov and Shekhtman [2] proved the inequality
$$
\max\limits_{\omega^N=1}|P(\omega)|\geq \left(1+C
\log\frac{N}{N-n}\right)^{-1} \max\limits_{|z|=1}|P(z)|,\eqno(2)
$$
where the absolute constant $C$ can be estimated by the number 16
(see [2, p. 3, 5]). This result generalizes an estimate due to
Marcinkiewicz [3], who obtained (2) for the case $N = n+ 1$.
Inequality (2) is better than inequality (1) for $n/N$ close to 1,
but worse than the Sheil-Small estimate for small values of $n/N$.
In the present paper, we prove the following statement.

\vskip0.5cm\textbf{Theorem.} \textit{Let $P$ be a polynomial of
degree $n$, and let $N$ be a natural number, $N\geq 2n$. Then, for
the discrete norm of the polynomial $P$, the following inequality
holds:}
$$
\max\limits_{\omega^N=1}|P(\omega)|\geq\cos\frac{\pi
n}{2N}\max\limits_{|z|=1}|P(z)|.
$$
\textit{The equality in (3) is attained in the case
$P(z)=(z\exp(i\pi/N))^n+1$ and for any $N$ which is a multiple of
$n$.}

\vskip0.5cm Estimate (3) holds for all $N >n$. However, for $n < N
< 2n$, it is worse than estimate (1). In the case $N = 2n$,
estimates (1) and (3) coincide and, for $N >2n$, inequality (3)
strengthens inequalities (1) and (2). Moreover, for numbers $N$
which are multiples of $n$, inequality (3) is sharp, and we obtain
the equality
$$
\sup\left\{\left(\max\limits_{|z|=1}|P(z)|\right)/\left(\max\limits_{\omega^N=1}
|P(\omega)|\right)\right\}=\left(\cos\left(\frac{\pi
n}{2N}\right)\right)^{-1},
$$
where $N$ is a multiple of $n$, and the upper bound is taken over
all polynomials $P$ of degree $n$ (see [2, Theorem 1]).

\vskip0.5cm Note that, as far back as 1931, Bernstein [4] obtained
inequalities for trigonometric sums close to inequalities (1)–(3)
(see also [5, pp. 147, 149, 154]). In particular, the corollary on
p. 154 of [5] implies inequality (3) in the case of even degrees
$n$. Our proof is different from the proofs in [1]–[5]. It is
based, essentially, only on the maximum principle of the modulus
of the suitable analytic function. Following [6], we can
strengthen inequality (3) by taking the constant term and the
leading coefficient of the polynomial $P$ into account.

\begin{center}2. AUXILIARY RESULT\end{center}

We introduce the notation
$$
m(P)=\min\limits_{|z|=1}|P(z)|,\qquad\quad
M(P)=\max\limits_{|z|=1}|P(z)|.
$$
We shall need the following analog of the Schwartz lemma in one
its particular cases. \vskip0.5cm\textbf{Lemma.} \textit{Let $P$
be a polynomial of degree $n$ for which $P(0)\neq 0$ and $m(P)\neq
M(P)$, and let the function $\zeta=\Phi(w)$ conformally and
univalently map the exterior of the closed interval
$\gamma=[m^2(P),M^2(P)]$ onto the disk $|\zeta|<1$ so that
$\Phi(\infty)=0$ and $\Phi(m^2(P))=-1$. Then the function}
$$
f(z)=\Phi\left( \overline{P(\overline{z})}P(\frac{1}{z})\right)
$$
\textit{is regular on the set}
$$
G=\left\{z:|z|<1,\;
\overline{P(\overline{z})}P(\frac{1}{z})\not\in\gamma\right\},
$$
\textit{analytically continuable to the set}
$$
E=\left\{z:|z|=1,\; |P(\overline{z})|\neq
m(P),\;|P(\overline{z})|\neq M(P)\right\},
$$
\textit{and, at the points of the set $E$, the following
inequality holds:}
$$
|f'(z)|\leq n.\eqno(4)
$$

\textbf{Proof.} The smoothness of the function $f$ on the sets $G$
and $E$ can easily be verified. Further, in a neighborhood of the
origin, the following expansion is valid:
$$
f(z)=\frac{M^2(P)-m^2(P)}{4\overline{c_0}c_n}z^n+\ldots\;,
$$
where $c_0$ is the constant term and $c_n$ is the leading
coefficient of the polynomial $P$. In addition, $f(z)\neq 0$ in
$G\setminus\{0\}$. Therefore, the function $z^n/f(z)$ is regular
on the open set $G$. At the points of the boundary of this set,
the modulus of this function does not exceed 1. By the maximum
principle for the modulus, we find that the inequality
$$
|f(z)|\geq|z^n|\eqno(5)
$$
holds everywhere on the set $G$. Now, let $z$ be an arbitrary
fixed point of the set $E$. Taking inequality (5) into account, we
obtain
$$
|f'(z)|=\frac{\partial|f(z)|}{\partial|z|}=\lim\limits_{r\to
1}\frac{|f(z)|-|f(rz)|}{1-r}\leq\lim\limits_{r \to
1}\frac{1-r^n}{1-r}=n.
$$
The lemma is proved.

\begin{center}3. PROOF OF THE THEOREM\end{center}

We can assume that $M(P)=1$ and $P(0)\neq 0$. Under these
conditions, $m(P)<M(P)$. Indeed, otherwise, the polynomial $P$
maps the circle $|z|=1$ into itself and, in view of the equality
$P(\infty)=\infty$, the symmetry principle leads to a
contradiction: $P(0)=0$. Let us show that, for any point
$z=e^{i\varphi}$ on the circle $|z|=1$, the following inequality
holds:
$$
\left|\left(|P(z)|^2\right)'_{\varphi}\right|\leq
n\sqrt{|P(z)|^2(1-|P(z)|^2)}\eqno(6)
$$
(see [6, Theorem 2]). In view of the continuity, it suffices to
verify this inequality for all points $z$ such that $\overline{z}$
belongs to the set $E$ from the lemma. Suppose that $z\in E$.
Then, for the function $f$ from the lemma, we have
$$
|f'(z)|=\left|\Phi'\left(\overline{P(\overline{z})}P(\frac1z)\right)\right|\left|\overline{P'
(\overline{z})}P(\frac1z)-\frac{1}{z^2}\overline{P(\overline{z})}P'(\frac1z)\right|=
$$
$$
=\left|\Phi'\left(|P(\overline{z})|^2\right)\right|\left|\frac{P^2(\overline{z})}{z}\right|\left|
\frac{\overline{\overline{z}P'(\overline{z})}}{\overline{P(\overline{z})}}-\frac{\dfrac1zP'(\dfrac1z)}{
P(\dfrac1z)}\right|=
$$
$$
=\left|\Phi'\left(|P(\overline{z})|^2\right)\right||P^2(\overline{z})|\left|
2{\rm
Im}\frac{\overline{z}P'(\overline{z})}{P(\overline{z})}\right|=
\left|\Phi'\left(|P(\overline{z})|^2\right)\right|\left|\left(|P(\overline{z})|^2\right)'_{\varphi}\right|.
$$
Before calculating the derivative $\Phi'$, note that the inverse
function $\Phi^{-1}(\zeta)$ is of the form
$$
\Phi^{-1}(\zeta)=\frac14\left(\zeta+\frac{1}{\zeta}\right)
(M^2(P)-m^2(P))+ \frac12(M^2(P)+m^2(P)).
$$
Hence
$$
\left|\Phi'\left(|P(\overline{z})|^2\right)\right|^{-1}=
\left|\frac14 \left(1-e^{-i2\theta}\right)(M^2(P)-m^2(P))\right|=
\frac12|\sin \theta|(M^2(P)-m^2(P)),
$$
where $\Phi^{-1}(e^{i\theta})=|P(\overline{z})|^2$ and, therefore,
$$
\cos
\theta=\frac{2|P(\overline{z})|^2-(M^2(P)+m^2(P))}{M^2(P)-m^2(P)}.
$$
Finally,
$$
|f'(z)|=\frac{\left|\left(|P(\overline{z})|^2\right)'_{\varphi}\right|}{
\sqrt{(|P(\overline{z})|^2-m^2(P))(M^2(P)-|P(\overline{z})|^2)}}.
$$
Using inequality (4), we obtain inequality (6) for $\overline{z}$.

Let us now pass to the proof of inequality (3). Let
$z_0=e^{i\varphi_0}$ denote one of the points at which the maximum
$M(P)=|P(z_0)|=1$ is attained, and let $\omega_k$ be the $N$th
root of 1 for which the arc of the circle
$$
\left\{z:|z|=1,\;|\arg z-\arg\omega_k|\leq\frac{\pi}{N}\right\}
$$
contains the point $z_0$. Suppose that, for some branch of the
argument, the following inequality holds:
$$
\arg\omega_k-\frac{\pi}{N}\leq\varphi_0\leq\arg \omega_k.
$$
Dividing both sides of inequality (6) by the quadratic root on the
right and integrating the resulting relation on the interval
$(\varphi_0,\arg\omega_k)$ with the replacement
$|P(z)|^2=u(\varphi)$, we obtain
$$
n(\arg\omega_k-\varphi_0)\geq\int\limits_{\varphi_0}^{\arg\omega_k}
\frac{-u'_{\varphi}d\varphi}{\sqrt{u(1-u)}}=-\int\limits_{1}^{u_k}
\frac{du}{\sqrt{u(1-u)}}=
$$
$$
=-2\int\limits_{1}^{\sqrt{u_k}}\frac{dt}{\sqrt{1-t^2}}=-2\arcsin\sqrt{u_k}+\pi,
$$
where $u_k=u(\arg\omega_k)$. Hence
$$
2\arcsin|P(\omega_k)|\geq\pi-n\frac{\pi}{N}>0,
$$
and
$$
|P(\omega_k)|\geq\sin\left(\frac{\pi}{2}-\frac{n\pi}{2N}\right)=\cos\frac{\pi
n}{2N}.
$$
Passing to the maximum, we obtain inequality (3). In the case
$\arg\omega_k\leq\varphi_0$, similar arguments yield the same
inequality.

If $P(z)=(z\exp(i\pi/N))^n+1$ and $N=nl$, where $l\geq 2$ is a
natural number, then
$$
\max\limits_{|z|=1}|P(z)|=2.
$$
On the other hand, direct calculations yield
$$
\max\limits_{\omega^N=1}|P(\omega)|=\max\limits_{0\leq k\leq
l-1}|P(\omega_k)|= \max\limits_{0\leq k\leq l-1}
2\left|\cos(\frac{\pi}{l}k+\frac{\pi}{2l})\right| =
2\cos\frac{\pi}{2l}.
$$
Thus, for the given polynomial $P$, we have the equality sign in
(3). The theorem is proved.

\begin{center}REFERENCES\end{center}
\begin{enumerate}
\item T. Sheil-Small, Bull. London Math. Soc. \textbf{40} (6), 956 (2008).
\item E. Rakhmanov and B. Shekhtman, J. Approx. Theory \textbf{139} (1-2), 2
(2006).
\item J. Marcinkiewicz, Acta Litt. Sci. Szeged \textbf{8}, 131 (1937).
\item S. N. Bernstein, Izv. Akad. Nauk SSSR OMEN \textbf{9}, 1151 (1931).
\item S. N. Bernstein, Collected Works, Vol. 2: The Constructive Theory
ofFunctions (Izdat. Akad. Nauk SSSR, Moscow, 1954) [in Russian].
\item V. N. Dubinin,Mat. Sb. \textbf{191} (12), 51 (2000) [Russian Acad. Sci.
Sb.Math. \textbf{191} (12), 1797 (2000)].
\end{enumerate}

\end{document}